\documentclass[12pt]{article}
\usepackage[latin2]{inputenc}
\usepackage{amsfonts, marvosym}
\usepackage{amsmath}
\usepackage{enumerate}
\usepackage{amssymb}
\usepackage{theorem}
\usepackage{authblk}
\usepackage{hyperref}
\usepackage{mathtools}

\newcommand{\R}{{\ensuremath{\mathbb{R}}}}
\newcommand{\N}{{\ensuremath{\mathbb{N}}}}

\renewcommand{\dj}{d\kern-0.4em\char"16\kern-0.1em}

\newcommand{\proof}{\noindent\textbf{Proof.}\ }
\newcommand{\proofof}{\noindent\textbf{Proof of}\ }

\newcommand{\qed}{\hfill\ensuremath{\Box}\\}

\newcommand{\Hj}{\hyperlink{(H1)}{\bf(H1)}}
\newcommand{\Hd}{\hyperlink{(H2)}{\bf(H2)}}

\newcommand{\Hpsi}{H^{\psi,1}}

\newcommand{\supp}{\operatorname{supp}}

\newtheorem{Thm}{Theorem}[section]
\newtheorem{Prop}[Thm]{Proposition}
\newtheorem{Cor}[Thm]{Corollary}
\newtheorem{Lem}[Thm]{Lemma}

\newtheorem{Def}[Thm]{Definition}

\newtheorem{Rem}[Thm]{Remark}
{\theorembodyfont{\upshape} }
\numberwithin{equation}{section}

\title{A note on the trace theorem for Besov-type spaces of generalized smoothness on $d$-sets\footnote{This work was supported by the Croatian Science Foundation under the project 3526.}}
\author{V.~Wagner\footnote{Department of Mathematics, University of Zagreb, 10000 Zagreb, Croatia\\\Letter\quad\href{mailto:wagner@math.hr}{wagner@math.hr}}}

\begin{document}

\maketitle

\begin{abstract}
The main goal of this paper is to give a complete proof of the trace theorem for Besov-type spaces of generalized smoothness associated with complete Bernstein functions satisfying certain scaling conditions on $d$-sets $D\subset\R^n$, $d\leq n$. The proof closely follows the classical approach by Jonsson, Wallin in \cite{jonwal} and the trace theorem for classical Besov spaces. Here, the trace space is defined by means of differences. When $d=n$, as an application of the trace theorem, we give a condition under which the test functions $C_c^\infty(D)$ are dense in the trace space on $D$.
\end{abstract}


\section{Introduction}

Besov-type spaces of generalized smoothness were introduced in the seventies by M.L. Goldman and G.A. Kalyabin as a generalization of the classical Sobolev and Besov spaces. Since then they have been studied by many authors from various points of view. Due to their close connection to the theory of stochastic processes and probability theory, these spaces continue to be of further interest. For a unified and general approach to function spaces of generalized smoothness and a review of results we refer the reader to \cite{farkas} and \cite{farkasleopold}.

In this paper we consider traces on $d$-sets of spaces of generalized smoothness associated with complete Bernstein functions and their representation via differences. This approach allows us to easily associate these spaces with domains of Dirichlet forms corresponding to a certain class of purely discontinuous symmetric Markov processes. A trace theorem on $h$-sets for generalized Besov-type spaces associated with functions satisfying slightly different assumptions was presented by Knopova and Z\"ahle in \cite{knop}. The authors consider the quarkonial representations of the trace spaces, as well as the representation via differences, and prove their equivalence. In the following paper \cite{zahle}, Z\"ahle proves a potential representation of the trace spaces and the corresponding trace Dirichlet forms. In order to make the theory of generalized Besov spaces more approachable in the context of application to the theory of stochastic processes, we present a complete and detailed proof of the trace theorem based on the characterization via differences, following the approach in \cite{jonwal}. 

Let $j:(0,\infty)\to(0,\infty)$ be a non-increasing function satisfying
\begin{equation}\label{eq:j1}
j(r) \leq c_1 j(r + 1) \text{ and }\int_0^\infty (1\wedge s^2)j(s)ds\leq\infty. 
\end{equation}
for all $r \geq 1$ and some constant $c_1>0$. Such a function $j$ is a radial L\'evy density of an isotropic unimodal L\' evy process with characteristic exponent
\begin{equation}\label{eq:j3}
 \psi(|\xi|)=\int_{\R^n\setminus\{0\}}\left(1-\cos{(x\cdot\xi)}\right)j(|x|)dx,\ \xi\in\R^n. 
\end{equation}  
Throughout the paper we will assume that $\psi$ satisfies one or both of the following scaling conditions, \\
{\hypertarget{(H1)}{\bf(H1)}}: There exist constants $0 < \delta_1 \leq \delta_2 < 1$ and $a_1,a_2 > 0$ such that
$$a_1 \lambda^{2\delta_1} \psi(t) \leq \psi(\lambda t) \leq a_2 \lambda^{2\delta_2} \psi(t),\quad \lambda\geq 1,\, t \geq1;$$
{\hypertarget{(H2)}{\bf(H2)}}: There exist constants $0 < \delta_3 \leq \delta_4 < 1$ and $a_3,a_4 > 0$ such that
$$a_3 \lambda^{2\delta_3} \psi(t) \leq \psi(\lambda t) \leq a_4 \lambda^{2\delta_4} \psi(t),\quad \lambda\geq 1,\, t < 1.$$
Under condition \Hj, by \cite[(2.1), (2.2)]{martin}, there exists a complete Bernstein function $\phi$ and a constant $\gamma_2\geq1$ such that
\begin{equation}\label{eq:phipsi}
 \gamma_2^{-1}\phi(|\xi|^2)\leq \psi(|\xi|)\leq \gamma_2\phi(|\xi|^2), \ \xi\in\R^n, 
\end{equation}
and the radial L\'evy density $j$ enjoys the following property: for every $R>0$
\begin{equation}\label{eq:j4}
 j(r)\asymp\frac{\phi(r^{-2})}{r^n},\ r\in(0,R).
\end{equation}
Due to the equivalence of norms, we will always assume that the constant $\gamma_2$ in \eqref{eq:phipsi} is equal to 1. Furthermore, by \cite[Lemma 2.1]{martin} every Bernstein function $\phi$ satisfies the following useful inequality,
\begin{equation}
  1\wedge\lambda\leq\frac{\phi(\lambda r)}{\phi(r)} \leq 1 \vee\lambda,\quad \lambda, r>0.\label{eq:phi}
\end{equation}
For further results on Bernstein functions, we refer the reader to \cite{bernstein}. Define by
\begin{equation}\label{eq:besselspace}
 \Hpsi(\R^n)=\{u\in L^2(\R^n): ||u||_{\psi,1}:=\int_{\R^n}(1+\psi(\xi))|\mathcal F u(\xi)|^2d\xi<\infty\}
\end{equation}
the $\psi$-Bessel-type potential space on $\R^n$. This space naturally arises as the domain of the Dirichlet form associated with the L\' evy process with the characteristic exponent $\psi$. In the first chapter, we give basic definitions and results regarding Besov-type spaces of generalized smoothness, connect them with $\psi$-Bessel-type potential spaces and provide several characterizations of the corresponding norms, based on differences. 

In the following two chapters, we give a complete proof of the trace theorem for spaces $H^{\psi,\alpha}(\R^n)$ (see Definition \ref{eq:convolution}) on $d$-sets, where the trace space is defined by means of differences. 
\begin{Def}\label{dmeasure}
Let $D$ be a non-empty Borel subset of $\R^n$ and $0<d\leq n$. A positive Borel measure $\mu$ on $D$ is called a $d$-measure if there exist positive constants $c_1$ and $c_2$ such that for all $x\in D$ and $r\in (0, 1]$,
 $$c_1 r^d\leq \mu(B(x,r))\leq c_2 r^d.$$
A non-empty Borel set $D$ is called a $d$-set if there exists a $d$-measure $\mu$ on $D$. 
\end{Def}
Note that by definition all $d$-measures on $D$ are equivalent to the restriction of the $d$-dimensional Hausdorff measure to $D$. For a $d$-set $D$ in $\R^n$ with $d$-measure $\mu$, we define the trace space as
\begin{align}
H^{\psi,1}(D,\mu)&=\{u\in L^2(D,\mu):||u||_{(1),D,\mu}<\infty\},\nonumber\\
||u||_{(1),D,\mu}&=||u||_{L^2(D,\mu)}+\left(\,\,\iint\limits_{|x-y|<1}|u(x)-u(y)|^2 \frac{\phi\left(|x-y|^{-2}\right)}{|x-y|^{2d-n}}\mu(dx)\mu(dy)\right)^\frac 1 2.\nonumber
\end{align}
The proof of the trace theorem is divided into three parts; the restriction theorem, the extension theorem for $d<n$ and finally, the extension theorem when $d=n$. We state the main result:
\begin{Thm}\label{thm:main}
Let $D$ be a $d$-set in $\R^n$ and $\psi$ a radial function defined by \eqref{eq:j3} such that \Hj and \Hd hold.
 \begin{enumerate}[(i)]
  \item There exists a continuous restriction operators $R:\Hpsi(\R^n)\to H^{\psi,1}(D,\mu)$.
  \item There exists a continuous extension operator $E:H^{\psi,1}(D,\mu)\to H^{\psi,1}(\R^{n})$, such that $R E u=u$, $\mu$-a.e. on $D$, for all $u\in H^{\psi,1}(D,\mu)$.
 \end{enumerate}
\end{Thm}
This result is a direct consequence of Theorem \ref{restriction_theorem}, Theorem \ref{extension_theorem1} and Theorem \ref{extension_theorem3} for spaces $H^{\psi,\alpha}(\R^n)$. As a consequence of Theorem \ref{thm:main} and \cite[Theorem 1.1]{wagner2}, we arrive to the following result on the correspondence of spaces $H^{\psi,1}(D)$ and $H^{\psi,1}_0(D)$, where the latter is the closure of $C_c^\infty(D)$ with respect to the norm $||\cdot||_{(1),D,\lambda_D}$. For an equivalent result in the case of the classical Bessel potential space see \cite[Corollary 2.8]{cen}.

\begin{Cor}\label{corollary}
Suppose that $D\subset\R^n$ is an open $n$-set, $\phi$ is a complete Bernstein function such that \Hj and \Hd hold and $\psi(\xi)=\phi(|\xi|^2)$.
\begin{enumerate}[(i)]
 \item Suppose that $2\delta_2 \leq n$ and that $\mathcal H_h (\partial D \cap K_m ) < \infty$ for an increasing sequence of Borel sets $K_m$ such that $\cup_{m\in\N}K_m \supset \partial D$, where $h(r) = r
^{n-2\delta_2}$ if $2\delta_2 < n$ and $h(r) = \max\{\log r , 0\}$ when $2\delta_2 = n=1$, then $H^{\psi,1}(D)=H^{\psi,1}_0(D)$.

\item If $2(\delta_1\wedge \delta_3)\geq n=1 $ or $\mathcal H^d(\partial D)>0$ for some $d>n-2\delta_1\geq 0$ then $H^{\psi,1}_0(D)\subsetneq H^{\psi,1}(D)$.
\end{enumerate}
\end{Cor} 

For easier notation, we write $f\asymp g$ if there exists a constant $c>1$ such that for all $x$, $c^{-1}g(x)\leq f(x)\leq c g(x)$. Equivalently, $f\lesssim g$ if there exists a constant $c>0$ such that for all $x$, $f(x)\leq c g(x)$. By $\lambda$ we denote the Lebesgue measure on $\R^n$.

\section{Besov-type spaces of generalized smoothness and equivalent norms}\label{section:besov}

The aforementioned $\psi$-Bessel-type potential space is a type of a much more general class of function spaces called Besov-type spaces of generalized smoothness. First we recall these spaces in their most general form.     
 
\begin{Def} A sequence $(\gamma_j)_{j\in\N_0}$ of positive real numbers is called
\begin{enumerate}[(i)]
\item \emph{almost increasing} if there exists $d_0>0$ such that $d_0\gamma_j\leq\gamma_k$, for all $j\leq k;$
\item \emph{strongly increasing} if it is almost increasing and in addition there exists a $\kappa_0\in\N$ such that $2\gamma_j\leq\gamma_k$, for $j\leq k-\kappa_0;$
\item of \emph{bounded growth} if there are positive constants $d_1$ and $J_0\in\N_0$ such that $\gamma_{j+1}\leq d_1\gamma_j$, for all $j\geq J_0;$
 \item an \emph{admissible sequence} if both $(\gamma_j)_{j\in\N_0}$ and $(\gamma_j^{-1})_{j\in\N_0}$ are of bounded growth and $J_0=0$, i.e.~there exist positive constants $d_0$ and $d_1$ such that $d_0\gamma_j\leq\gamma_{j+1}\leq d_1\gamma_j$, for all $j\in\N_0.$
\end{enumerate}
\end{Def}
\begin{Def}
Let $N=(N_j)_{j\in\N_0}$ be a strongly increasing sequence and define
\begin{align*}
&\Omega_0^N=\{x\in\R^n:|x|\leq N_0\} \text{ and }\Omega_j^N=\{x\in\R^n:N_{j-1}\leq|x|\leq N_{j+1}\},\quad j\in\N.
\end{align*}
Let $\Phi^N$ be a collection of all partitions of unity of $C_c^\infty(\R^n)$ functions associated with this decomposition. Let $\sigma=(\sigma_j)_{j\in\N_0}$ be an admissible sequence respectively and $(\varphi_j^N)_{j\in\N_0}\in\Phi^N$. The Besov space of generalized smoothness associated with $N$ and $\sigma$ is defined by
$$B_2^{\sigma,N}=\{g\in S'(\R^n):||g||_{B,\sigma,N}:=||(\sigma_j\varphi^N_j(D)g)_{j\in\N_0}|l_2(L_2(\R^n))||<\infty\},$$
where $\varphi(D)g(x)=\mathcal F^{-1}(\varphi(\cdot)\mathcal F g)(x)$ and
$||(f_j)_{j\in\N_0}|l^2(L^2(\R^n))||=\left(\sum_{j=0}^\infty ||f_j||_{L^2(\R^n)}^2\right)^{\frac 1 2}.$
\end{Def}
By \cite[Remark 10.1.2.]{farkas} the space $B_2^{\sigma,N}$ is independent of the choice of system $(\varphi_j^N)_{j\in\N_0}$ in the sense of equivalent norms. We will restrict ourselves to a special subclass of spaces $B_2^{\sigma,N}$ associated with an admissible symbol.
\begin{Def}\label{def:admissible}
A non-negative function $a\in C^\infty(\R^n)$ is an \emph{admissible symbol} if
\begin{enumerate}[(i)]
\item $\lim\limits_{|x|\to\infty}a(x)=\infty$,
\item $a$ is almost increasing in $|x|$, i.e.~there exist constants $\delta_0\geq 1$ and $R>0$ such that $a(x)\leq\delta_0 a(y)$ if $R\leq |x|\leq |y|$,
\item there exists an $m>0$ such that $x\to \frac{a(x)}{|x|^m}$ is almost decreasing in $|x|$,
\item for every multi-index $\alpha\in\N_0^n$ there exist constants $c_\alpha>0$ and $R>0$ such that
$$|D^\alpha a(x)|\leq c_\alpha\frac{a(x)}{(1+|x|^2)^{|\alpha|/2}}, \quad \forall |x|\geq R.$$
\end{enumerate}
The family of all admissible functions will be denoted by $\mathcal A$.
\end{Def}
Note that for every Bernstein function $\phi$ such that $\lim\limits_{r\to\infty}\phi(r)=\infty$, the function $\xi\to\phi(|\xi|^2)$ is an admissible symbol, see \cite[Lemma 3.1.13]{farkasleopold}.
\begin{Rem}\label{trace:admissible}
 By \cite[Lemma 3.1.17, Remark 3.1.18]{farkasleopold}, the sequence $N_j^{a,r}=\sup\{|x|:a(x)\leq 2^{jr}\}$, $j\in\N_0,$ where $a\in\mathcal A$ and $r>0$, is strongly increasing.
\end{Rem}

For $a\in\mathcal A$ we define the Besov space of generalized smoothness associated with $a$ as $B^{\sigma,N^{a,2}}_2(\R^n),$ where $\sigma=\{2^{j}\}_{j\in\N_0}$ is an admissible sequence. These spaces have two useful representations in the sense of equivalent norms; one given by the Littlewood-Paley-type theorem and the other by means of differences.  

\begin{Prop}\label{trace:Ka}\cite[Theorem 3.1.20, Corollary 3.1.21]{farkasleopold}\\
Let $a\in\mathcal A$, $N=N^{a,2}$ the strongly increasing sequence associated with $a$, $\alpha>0$ and $\sigma^\alpha=\{2^{\alpha j}\}_{j\in\N_0}$ an admissible sequence. Then the space $(H^{a,\alpha}(\R^n),||\cdot||_{a,\alpha})$, defined by
\begin{align}\label{eq:Ka}
 &H^{a,\alpha}(\R^n)=\{u\in L^2(\R^n):||u||_{a,\alpha}<\infty\},\nonumber\\
 &||u||_{a,\alpha}:=||(id+a(D))^{\alpha/2}u||_{L_2(\R^n)}=\left(\int_{\R^n} (1+a(\xi))^\alpha|\hat u(\xi)|^2d\xi\right)^\frac 1 2,
\end{align}
is equivalent to $(B^{\sigma^\alpha,N^{a,2}}_2(\R^n),||\cdot||_{B,\sigma^\alpha,N^{a,2}})$.
\end{Prop}

Proposition \ref{trace:Ka} implies that Besov-type spaces of generalized smoothness associated with $\psi$ can be characterized as
\begin{align*}
& H^{\psi,\alpha}(\R^n)=\left\{u\in S'(\R^n):\exists f\in L^2(\R^n)\text{ such that }\mathcal F u=(1+\psi)^{-\alpha/2}\mathcal F f\right\}
\end{align*}
 Since the function $x\mapsto(1+x)^{-\alpha/2}$ is completely monotone for every $\alpha>0$, by \cite[Theorem 3.7]{bernstein} the function $(1+\phi)^{-\alpha/2}$ is also completely monotone. By Schoenberg's theorem, \cite[Theorem 2]{schoenberg} it follows that the function {${(1+\psi)^{-\alpha/2}}$} is a positive definite function and therefore a Fourier transforms of an integrable function, \cite[Theorem 4.14]{bernstein}, called the Bessel-type potential $K_{\psi,\alpha}$. This means that the space $H^{\psi,\alpha}(\R^n)$ defined by \eqref{eq:Ka} can be characterized as a convolution space via the $\psi$-Bessel convolution kernel $K_{\psi,\alpha}$ i.e.
\begin{equation}\label{eq:convolution}
H^{\psi,\alpha}(\R^n)=\{K_{\psi,\alpha}*f:f\in L^2(\R^n)\},\quad ||K_{\psi,\alpha}*f||_{\psi,\alpha}:=||f||_{L^2(\R^n)}.
\end{equation}
From now on we assume that the function $\psi$ satisfies conditions \Hj and \Hd. The following estimates for the kernels $K_{\psi,\alpha}$ were obtained in \cite[Remark 33, Remark 34]{knop} and \cite{zahle}, for further results on Bessel-type potential spaces we refer the reader to the latter. 
\begin{Lem}\label{trace:Kpsi}
Let $\alpha>0$ and assume \Hj and \Hd hold. If $\alpha(\delta_2\vee\delta_4)<n$ then there exist constants $c_i=c_i(\phi,\alpha, n)>0$, $i=1,2$, such that for all $x\in \R^n$ and $0\leq j\leq n$
\begin{align}
|K_{\psi,\alpha}(x)|&\leq \frac{c_1}{|x|^n\phi(|x|^{-2})^{\alpha/2}},\nonumber\\
|(K_{\psi,\alpha}(x))'_{x_j}|&\leq \frac{c_{2}}{|x|^{n+1}\phi(|x|^{-2})^{\alpha/2}}.\nonumber
\end{align}
\end{Lem}

Next, we consider the characterization of spaces $H^{\psi,\alpha}(\R^n)$ via differences. First note that by applying \eqref{eq:j3} and the Parseval's identity (see for example \cite[Example 1.4.1]{fukushima}), it easily follows that 
\begin{align}\label{eq:equivalence1}
||u||_{\psi,1}\asymp ||u||_{L^2(\R^n)}+\left(\frac 1 2\int_{\R^n}\int_{\R^n} (u(x+y)-u(x))^2 j(y)dy dx\right)^\frac 1 2
\end{align}
for all $u\in\Hpsi(\R^n)$. We also introduce an equivalent norm on $\Hpsi(\R^n)$, which we later use in the proof of the trace theorem,
$$||u||_{(1)}=||u||_{L^2(\R^n)}+\left(\,\,\iint\limits_{|x-y|<1}|u(x)-u(y)|^2 \frac{\phi\left(|x-y|^{-2}\right)}{|x-y|^n} dx dy\right)^\frac 1 2,$$
with the equivalence of norms $||\cdot||_{(1)}$ and $||\cdot||_{\psi,1}$ following from \eqref{eq:equivalence1}, \eqref{eq:j4} and the fact that $j$ is a L\'evy measure. In the remainder of this chapter we prove an equivalent result for spaces $H^{\psi,\alpha}(\R^n)$.

\begin{Def}
For a function $f$ on $\R^n$, $h\in\R^n$ and $k\in\mathbb N$, the $k$-th difference of function $f$ is defined by $(\Delta_h^kf)(x):=\Delta_h^1(\Delta_h^{k-1}f)(x)$, $x\in\R^n$, where $\Delta_h^1f(x)=f(x+h)-f(x)$. Next, define the $k$-th modulus of continuity of a function $f\in L^2(\R^n)$ as $\omega_k(f,t)=\sup_{|h|<t}||\Delta_h^kf||_{L^2(\R^n)}$, $t>0$. Furthermore, for an admissible sequence $(\gamma_j)_{j\in\N_0}$ let $\overline\gamma_j=\sup_{k}\frac{\gamma_{j+k}}{\gamma_k}$ and $\underline\gamma_j=\inf_{k}\frac{\gamma_{j+k}}{\gamma_k}$, and let $\underline{s}(\gamma):=\lim\limits_{j\to\infty}j^{-1}\log\underline\gamma_j$ and $\overline{s}(\gamma):=\lim\limits_{j\to\infty}j^{-1}\log\overline\gamma_j$ be the lower and upper Boyd index, respectively.
\end{Def}
Since $\gamma_{j+i+k}\leq \overline\gamma_j\gamma_{i+k}$ for all $i,j,k\in\N_0$ it follows that $\overline\gamma_{j+i}\leq \overline\gamma_j\overline\gamma_{i}$, so the sequence $\log\overline\gamma_j$ is subadditive. By Fekete's subadditive lemma the sequence $\left(\frac{\log\overline\gamma_j}{j}\right)_j$ converges to $\inf\limits_j\frac{\log\overline\gamma_j}{j}$, so the upper index $\overline{s}(\gamma)$ is well defined. The analogous conclusion follows for the lower index $\underline{s}(\gamma)$, since $\log\underline\gamma_j=-\log\left(\overline{\gamma^{-1}}_j\right)$. 

\begin{Thm}\cite[Theorem 4.1]{moura}\label{trace:moura}\\
Let $\sigma$ and $N$ be admissible sequences and $\underline{N}_1=\inf\limits_{k\geq 0}\frac{N_{k+1}}{N_k}>1$ and $\displaystyle{\frac{\underline{s}(\sigma)}{\overline{s}(N)}}>0$. Let $k$ be an integer such that $\displaystyle{k>\frac{\overline{s}(\sigma)}{\underline{s}(N)}}$. Then the norm $||\cdot||_{B,\sigma,N}$ on $B_2^{\sigma,N}$ is equivalent to 
\[
||u||_{L^2(\R^n)}+\left(\sum_{j=0}^\infty\sigma_j^2\omega_k(u,N_j^{-1})^2\right)^\frac{1}{2}.
\]
\end{Thm}

For a similar result see also \cite[Theorem 16]{knop}. Note that for every $a\in\mathcal A$ the sequence $N^{a,r}$ satisfies the assumption $\underline{N}^{a,r}_1>1$. Furthermore, recall that $\sigma^\alpha=(2^{\alpha j})_{j\in\N_0}$ is an admissible sequence and by Remark \ref{trace:admissible} the sequence $N^{\psi,2}$ is strongly increasing. One can easily show that $N^{\psi,2}$ is also admissible. Furthermore,
\begin{align*}
	\frac{\underline{s}(\sigma^\alpha)}{\overline{s}(N^{a,2})}\geq\frac{\alpha}{2}>0 \mbox{  and  }	\frac{\overline{s}(\sigma)}{\underline{s}(N^{\psi,2})}\leq\frac{\alpha\log 2}{\frac{1}{\delta_2}\log 2}=\alpha\delta_2
\end{align*}
so Theorem \ref{trace:moura} holds for $k>\alpha\delta_2$. Furthermore,
\begin{align}
\sum_{j=0}^\infty 2^{\alpha j}\sup_{|H|<1/N_j^{\psi,2}}||\Delta_H^ku||_{L^2(\R^{n})}^2&\asymp \sum_{j=0}^\infty\int\limits_{2^{-(j+1)}\leq t<2^{-j}}\frac{1}{t^{1+2\alpha}}\sup_{|H|<1/N_j^{\psi,2}}||\Delta_H^ku||_{L^2(\R^{n})}^2dt\nonumber\\
&\asymp  \int_0^1\frac{1}{t^{1+2\alpha}}\sup_{|H|<1/\psi^{-1}(t^{-2})}||\Delta_H^ku||_{L^2(\R^{n})}^2dt,\label{eq:sup_H}
\end{align}
since by $2^{-(j+1)}\leq t<2^{-j}$ implies $\psi^{-1}(t^{-2})\asymp N_j^{\psi,2}.$ 
By change of variable $t^{-2}=\psi(|h|^{-1})$ it follows that \eqref{eq:sup_H} is comparable to
\[
\int\limits_{|h|<1}\frac{(\psi^\alpha)'(|h|^{-1})}{|h|^{n+1}}\sup_{|H|<|h|}||\Delta_H^ku||_{L^2(\R^{n})}^2dh.
\]
Since $(\psi^\alpha)'(t)\asymp\frac{\psi^\alpha(t)}{t}$ it follows that the last line is comparable to
\begin{align}
&\int\limits_{|h|<1}\frac{\psi^\alpha(|h|^{-1})}{|h|^{n}}\sup_{|H|<|h|}||\Delta_H^ku||_{L^2(\R^{n})}^2dh.\label{trace:a_norm_sup}
\end{align}
\begin{Rem}\label{trace:Delta^k}
\begin{enumerate}[(i)]
\item By applying a straightforward generalization of \cite[Theorem 2.6.1]{triebel2} to \eqref{trace:a_norm_sup}, Theorem \ref{trace:moura}, Theorem \ref{trace:Ka} and the calculation above imply that the norms
\begin{align}\label{trace:a_norm}
&||u||_{(1),\alpha,k}:=||u||_{L^2(\R^{n})}+\left(\int\limits_{|h|<1}\frac{\psi^\alpha(|h|^{-1})}{|h|^{n}}||\Delta_h^k u||_{L^2(\R^{n})}^2dh\right)^\frac{1}{2}
\end{align}
are equivalent to $||\cdot||_{a,1}$, for all $k> \alpha\delta_2$.
\item Since the function $\frac{\psi^\alpha(|\,\cdot\,|^{-1})}{|\,\cdot\,|^{n}}$ is continuous and $||\Delta_h^k u||_{L^2(\R^{n})}\leq c(k)|| u||_{L^2(\R^{n})}$ the norms ${||\cdot||_{(1),\alpha,k}^{h_0}}$, 
\[
||u||_{(1),\alpha,k}^{h_0}:=||u||_{L^2(\R^{n})}+\left(\int\limits_{|h|<h_0}\frac{\psi^\alpha(|h|^{-1})}{|h|^{n}}||\Delta_h^k u||_{L^2(\R^{n})}^2dh\right)^\frac{1}{2},
\]
are equivalent for all $h_0>0$.
\item Let $c>0$, $N\in\mathbb Z$ and $k> \alpha\delta_2$. The norm $||\cdot||_{(2)}^{c,N}$ on $H^{\psi,\alpha}(\R^n)$
defined by 
\begin{equation}\label{eq:norm2}
||u||_{(2),\alpha,k}^{c,N}=||u||_{L^2(\R^n)}+\left(\,\,\sum_{j=N}^\infty
 \psi^\alpha\left(2^{j}\right) 2^{nj}
 \int\limits_{|h|<c2^{-j}}||\Delta_h^k u||_{L^2(\R^{n})}^2dh\right)^\frac 1 2
\end{equation}
is equivalent to the norm $||\cdot||_{(1)}$. This follows by applying \eqref{eq:phi} to the norm in (ii) for $h_0=c2^{-N}$,
\begin{align}
&\int\limits_{|h|<c2^{-N}}\frac{\psi^\alpha(|h|^{-1})}{|h|^{n}}||\Delta_h^k u||_{L^2(\R^{n})}^2dh=\sum_{j=N}^\infty\,\,\int\limits_{c2^{-j-1}\leq|h|<c2^{-j}}\frac{\psi^\alpha(|h|^{-1})}{|h|^{n}}||\Delta_h^k u||_{L^2(\R^{n})}^2dh\nonumber\\
&\asymp\sum_{j=N}^\infty
\psi^\alpha\left(2^{j}\right)2^{nj}\int\limits_{c2^{-j-1}\leq|h|<c2^{-j}}||\Delta_h^k u||_{L^2(\R^{n})}^2dh\asymp \sum_{j=N}^\infty \sum_{i=N}^j \psi^\alpha\left(2^{i}\right) 2^{ni}
 \int\limits_{c2^{-j-1}\leq|h|<c2^{-j}}||\Delta_h^k u||_{L^2(\R^{n})}^2dh\nonumber\\
 &=  \sum_{i=N}^\infty \psi^\alpha\left(2^{i}\right) 2^{ni}
 \int\limits_{|h|<c2^{-j}}||\Delta_h^k u||_{L^2(\R^{n})}^2dh\nonumber.
\end{align}
\end{enumerate}  
\end{Rem}

\section{The restriction theorem} \label{section:restriction}

In this section we provide a detailed proof of the continuity of the restriction operator, as a generalization of \cite[Section V.1.2]{jonwal}. The same approach is used in \cite[Appendix III]{knop}. Before we start with the proof, we show the following useful consequence of the estimates on Bessel-type potentials from Lemma \ref{trace:Kpsi}. 

\begin{Lem}\label{trace:R1}
Let $d\leq n$, $D$ a $d$-set in $\R^n$ and $\mu$ the $d$-measure on $D$. Let $\psi$ be a function such that \eqref{eq:j3}, \Hj and \Hd hold and $\alpha>0$ such that
\begin{equation}\label{eq:conditionR1}
\frac{n-d}2<\alpha\delta_1\leq\alpha(\delta_2\vee\delta_4)<\frac{n-d}{2}+1.
\end{equation}
Then there exists a constant $c>0$ such that for all $r\leq \frac 1 3$ and $f\in L^2(\R^n)$
\begin{equation*}
\iint_{|x-y|<r} (K_{\psi,\alpha}\ast f(x)-K_{\psi,\alpha}\ast f(y))^2 \mu(dx) \mu(dy) \leq c \frac{r^{2d-n}}{\phi^\alpha(r^{-2})} ||f||_{L^2(\R^n)}^2.
\end{equation*}
\end{Lem}
\proof
Note that for every constant $0<a<1$
\begin{align}
(K_{\psi,\alpha}*f(x)-K_{\psi,\alpha}*f(y))^2\leq\int
|K_{\psi,\alpha}(x-t)-K_{\psi,\alpha}(y-t)|^{2a}f^2(t)dt\nonumber\cdot\int
|K_{\psi,\alpha}(x-t)-K_{\psi,\alpha}(y-t)|^{2(1-a)}dt. \nonumber
\end{align}
Let $|x-y|<r$. By Lemma \ref{trace:Kpsi} it follows that
\begin{align}
&\int_{|y-t|<2r}|K_{\psi,\alpha}(x-t)-K_{\psi,\alpha}(y-t)|^{2(1-a)}
dt\lesssim \int_{|z|<3r}|K_{\psi,\alpha}(z)|^{2(1-a)}dz\lesssim \int_{|z|<3r}\left(\frac 1{|z|^n\phi^{\frac\alpha 2}(|z|^{-2})}\right)^{2(1-a)}dz\nonumber\\
&\overset{\text{\Hj}}{\lesssim} \frac{r^{-2\alpha\delta_1(1-a)}}{(\phi^\frac\alpha 2(r^{-2}))^{2(1-a)}}\int_{|z|<3r}\left(\frac 1{|z|^{n-\alpha\delta_1}}\right)^{2(1-a)}dz = \tilde c_1 \frac{r^n}{\left(r^n\phi^{\frac \alpha 2}(r^{-2})\right)^{2(1-a)}} \nonumber,
\end{align}
for some $\tilde c_>0$ and $a$ such that
\begin{align}
&2(1-a)(n-\alpha\delta_1)<n\label{trace:a1}.
\end{align}
Analogously, if
\begin{align}
&2a(n-\alpha\delta_1)<d\label{trace:a2}
\end{align}
then for all $t\in\R^n$ there exists a constant $\tilde c_2>0$
\begin{align}
&\iint\limits_{\substack{|x-y|<r\\|y-t|<2r}}
|K_{\psi,\alpha}(x-t)-K_{\psi,\alpha}(y-t)|^{2a}\mu(dx)\mu(dy)\leq \tilde c_2\frac{r^{2d}}{\left(r^n\phi^{\frac\alpha 2}(r^{-2})\right)^{2a}}.\nonumber
\end{align}
Therefore, it follows that
\[
 \iint_{|x-y|<r}\left(\int_{|y-t|<2r} (K_{\psi,\alpha}(x-t)-K_{\psi,\alpha}(y-t))f(t)dt\right)^2\mu(dx)\mu(dy)\leq \tilde c_1\tilde c_2 \frac{r^{2d-n}}{\phi^\alpha(r^{-2})}||f||_{L^2(\R^n)}^2.
\]
For the second part, note that the mean value theorem and Lemma \ref{trace:Kpsi} imply that
\begin{align}\label{eq:Kpsilemma}
\iint\limits_{\substack{|x-y|<r\\|y-t|>2r}}
|K_{\psi,\alpha}(x-t)-K_{\psi,\alpha}(y-t)|^{2a}\mu(dx)\mu(dy)\lesssim r^{2a}\,\iint\limits_{\substack{|x-y|<r\\|y-t|>2r}}\left(\frac 1{|z_{x,y}|^{n+1}\phi^{\frac\alpha 2}(|z_{x,y}|^{-2})}\right)^{2a}\mu(dx)\mu(dy),
\end{align}
where $z_{x,y}=y-t+\theta_{x,y}(x-y)$ for some $\theta_{x,y}\in(0,1)$ and $|z_{x,y}|\geq \frac{|y-t|} 2$. Let $\delta=\delta_2\vee\delta_4$. By \Hj and \Hd there exists a constant $\tilde c_3>0$ such that the last line in \eqref{eq:Kpsilemma} is comparably less then
 \begin{align*}
\frac{r^{2a+d}}{\left(r^{\alpha\delta}\phi^{\frac \alpha 2}(r^{-2})\right)^{2a}}\int\limits_{|z|>2r}\left(\frac 1{|z|^{n+1-\alpha\delta }}\right)^{2a}\mu(dz)= \tilde c_{3}\frac{r^{2d}}{\left(r^n\phi^{\frac \alpha 2}(r^{-2})\right)^{2a}}
\end{align*}
if
\begin{equation}
2a(n+1-\alpha\delta )>d.\label{trace:a3}
\end{equation}
Similarly, if
\begin{equation}
2(1-a)(n+1-\alpha\delta )>n\label{trace:a4}
\end{equation}
then there exists a constant $\tilde c_{4}>0$ such that $|x-y|<r$ implies
\begin{equation}
\int\limits_{2r<|y-t|}|K_{\psi,\alpha}(y-t)-K_{\psi,\alpha}(x-t)|^{2(1-a)}dt\leq \tilde c_{4}\frac{r^n}{\left(r^n\phi^{\frac\alpha 2}(r^{-2})\right)^{2(1-a)}}.\nonumber
\end{equation}
Therefore, it follows that
\[
 \iint_{|x-y|<r}\left(\int_{|y-t|>2r} (K_{\psi,\alpha}(x-t)-K_{\psi,\alpha}(y-t))f(t)dt\right)^2\mu(dx)\mu(dy)\leq \tilde c_3\tilde c_4 \frac{r^{2d-n}}{\phi^\alpha(r^{-2})}||f||_{L^2(\R^n)}^2.
\]
Since $\frac{n-d}{2}<\alpha\delta_1\leq\alpha\delta<\frac{n-d}{2}+1$ it follows that 
$$\left(\frac{d}{2(n+1-\alpha\delta)}, 1-\frac{n}{2(n+1-\alpha\delta)}\right)\cap\left(1-\frac{n}{2(n-\alpha\delta_1)},\frac{d}{2(n-\alpha\delta_1)}\right)\neq\emptyset,$$ 
so we can choose $a$ such that \eqref{trace:a1}, \eqref{trace:a2}, \eqref{trace:a3} and \eqref{trace:a4} hold.
\qed

For a $d$-set $D$ in $\R^n$ with $d$-measure $\mu$  and $\alpha>0$ we define the trace space in terms of differences
\begin{align}
&H^{\psi,\alpha}(D,\mu)=\{u\in L^2(D,\mu):||u||_{(1),D,\mu}<\infty\},\nonumber\\
&||u||_{(1),D,\mu}=||u||_{L^2(D,\mu)}+\left(\,\,\iint\limits_{|x-y|<1}|u(x)-u(y)|^2 \frac{\psi^\alpha\left(|x-y|^{-1}\right)}{|x-y|^{2d-n}}\mu(dx)\mu(dy)\right)^\frac 1 2.\nonumber
\end{align}
Note that, similarly as in Remark \ref{trace:Delta^k}(iii), one can easily show that
\begin{equation*}
||u||_{(1),D,\mu}\asymp ||u||_{(2),D,\mu}^{c,N}=||u||_{L^2(D,\mu)}+\left(\,\,\sum_{j=N}^\infty
 \psi^\alpha\left(2^{j}\right) 2^{(2d-n)j}
 \iint\limits_{|x-y|<c2^{-j}}|u(x)-u(y)|^2
 \mu(dx)\mu(dy)\right)^\frac 1 2.
\end{equation*}
Next we define the restriction operator. Let $Ru$ be the pointwise restriction on $D$ of the {strictly defined} function corresponding to ${u\in H^{\psi,\alpha}(\R^{n})}$, i.e. 
$$Ru(x)=\lim_{r\to0}\frac 1 {\lambda(B(x,r))}\int_{B(x,r)}u(y)dy,\, x\in D,$$
whenever the limit exists. Here, instead of $u$ one can a take a quasi continuous modification $\tilde u$ of $u$, see for example \cite[Theorem 3.5.7]{farkas}. Recall that a function $\tilde u$ is quasi continuous if for every $\varepsilon>0$ there exists an open set $G_\varepsilon$ such that $\tilde u_{|G_\varepsilon^c}$ is continuous, $\text{Cap}_{\psi,\alpha}(G_\varepsilon)<\varepsilon$ and $\tilde u = u$ a.e. This means that outside of set $N=\cap_{\varepsilon>0} G_\varepsilon$ function $u$ can be strictly defined and therefore $\tilde Ru=u$ on $N^c$. To show that $\mu(N)=0$, note that by the definition of capacity, for every $\varepsilon>0$ there exists a function $v\in H^{\psi,\alpha}(\R^n)$ such that $v\geq 1$ on $N$ and $||v||_{a,1}<\varepsilon$, so
\[
\mu(N)\leq \left(\int_N |v(x)|^2\mu(dy)\right)^\frac 1 2\leq  \tilde c_1 ||v||_{\psi,\alpha}<\tilde c_1\varepsilon
\]
and therefore $\mu(N)=0$. Here the second inequality follows from calculations analogous to ones in the proof of Lemma \ref{trace:R1}.

\begin{Thm}\label{restriction_theorem}
Let $d\leq n$, $D$ a $d$-set in $\R^n$ and $\mu$ the $d$-measure on $D$. Let $\psi$ be a function such that \eqref{eq:j3}, \Hj and \Hd hold and $\alpha>0$ satisfying \eqref{eq:conditionR1}. There exist a continuous restriction operators $R:H^{\psi,\alpha}(\R^n)\to H^{\psi,\alpha}(D,\mu)$.
\end{Thm}
\proof  Using the classical interpolation theorem for a special class of spaces associated with $\Hpsi(\R^n)$ we will show that there exists a constant $\tilde c_1>0$ such that for all $u\in H^{\psi,\alpha}(\R^n)$
 \begin{equation}
  ||Ru||_{(2),D,\mu}^{1,1}\leq \tilde c_1 ||u||_{\psi,\alpha}.\label{trace:show2}
 \end{equation}
Denote by $a_j(x,y)=|Ru(x)-Ru(y)|1_{|x-y|<2^{-j}}$, $j\in \N_0$ and let $T$ be an operator on $L^2(\R^n)$ such that $Tu=(a_j)_{j\in\N_0}$. Since $Ru=u$ $\mu$-a.e. on $D$, by Lemma \ref{trace:R1} there exists a constant $\tilde c_{2}>0$ such that
\begin{align*}
\sup_{j\in\N_0} \left(\phi^\alpha\left(2^{2j}\right) 2^{j(2d-n)}
\iint\limits_{|x-y|<2^{-j}}|Ru(x)-Ru(y)|^2\mu(dx)\mu(dy)\right)\leq \tilde c_{2} ||f||_{L^2(\R^n)}^2
\end{align*}
for all $u=K_{\psi,\alpha}*f$, $f\in L^2(\R^n)$. Let $L=L^2(D\times D,\mu\times\mu)$. It follows that 
$$(a_j)_j\in l_\infty^{\phi,\alpha}(L)=\{(\xi_j)_{j\in\N_0}: \xi_j\in L,||\xi||_{l_\infty^{\phi,\alpha}(L)}=\sup_{j\in\N_0}\phi^\alpha\left(2^{2j}\right) 2^{nj}||\xi_j||_{L}^2<\infty\}$$
and that the operator $T$ is bounded from $H^{\psi,\alpha}(\R^n)$ to $l_\infty^{\phi,\alpha}(L)$, i.e.  
\[
||(a_j)_{j\in\N_0}||_{l_\infty^{\phi,\alpha}(L)}\leq \tilde c_{2}||K_{\psi,\alpha}*f||_{\psi,\alpha}=\tilde c_{2}||f||_{L^2(\R^n)}.
\]
Let $(X_1,X_2)_{\theta,2}:=\{a: \,a\in X_1+X_2,\, ||a||_{(X_1,X_2)_{\theta,2}} <\infty\}$ be the $K$-interpolation space of Banach spaces $X_1$ and $X_2$ and $||a||_{(X_1,X_2)_{\theta,2}}=\left(\int_0^\infty (t^{-\theta}K(t,a))^2\frac{dt}{t}\right)^\frac 1 2$, where the Peetre $K$-functional is defined by $K(t,a)=\inf\limits_{a=a_1+a_2}(||a_1||_{X_1}+t||a_2||_{X_2})$, see \cite[Section 1.3, p.23]{triebel}. Choose $\alpha_0<\alpha<\alpha_1$ such that $\frac{n-d}{2}<\delta_1\alpha_0\leq(\delta_2\vee\delta_4)\alpha_1<\frac{n-d}{2}+1$ and let $\theta=\frac{\alpha_1-\alpha}{\alpha_1-\alpha_0}\in(0,1)$. By \cite[Theorem 1.3.3(a)]{triebel} and \cite[Lemma 4.1]{interpolation}, since $T$ is bounded from $H^{\psi,\alpha_i}(\R^n)$ to $l^{\phi,\alpha_i}_\infty(L)$, $i=0,1$ it is also bounded from $(H^{\psi,\alpha_0}(\R^n),H^{\psi,\alpha_1}(\R^n))_{\theta,2}$ to $(l^{\phi,\alpha_0}_\infty(L),l^{\phi,\alpha_1}_\infty(L))_{\theta,2}$. By a version of \cite[Theorem 1.18.2]{triebel}, $(l_\infty^{\phi,\alpha_0}(L),l_\infty^{\phi,\alpha_1}(L))_{\theta,2}=l_2^{\phi,\alpha}(L)$, where
\[
l_2^{\phi,\alpha}(L)=\left\{(\xi_j)_{j\in\N_0}: \xi_j\in L,||\xi||_{l_2^{\phi}(L)}=\left(\sum_{j\in\N_0}\phi^\alpha\left(2^{2j}\right) 2^{nj}||\xi_j||_{L}^2\right)^\frac 1 2<\infty\right\}.
\]
Furthermore, $(l_2^{\phi,\alpha_0}(L^2(\R^n)),l_2^{\phi,\alpha_1}(L^2(\R^n)))_{\theta,2}=l_2^{\phi,\alpha}(L^2(\R^n))$ and $H^{\psi,\alpha}(\R^n)$ is a retract of the space $l_2^{\phi,\alpha}(L^2(\R^n))$, \cite[Theorem 2.5 and Theorem 3.4]{interpolationCobos}. Therefore, the interpolation identity
$(H^{\psi,\alpha_0}(\R^n),H^{\psi,\alpha_1}(\R^n))_{\theta,2}=H^{\psi,\alpha}(\R^n)$ follows by \cite[Theorem 5.3]{interpolationCobos}. Combining these results we get \eqref{trace:show2}.
 
\qed

\section{Extension theorem}

The proof of the extension theorem is divided into two parts; case $d<n$ and case $d=n$. In the first case we define the extension operator by using the Whitney decomposition of $D^c$ and the approach as in \cite[V.1.3]{jonwal}, which deals with the classical Besov spaces. The case $d=n$ follows then by interpreting the $n$-set $D$ in $\R^n$ as a $n$-subset of $\R^{n+1}$. The extension operator is then defined as a composition of the extension operator to $\R^{n+1}$ and restriction operator to $\R^{n}$. 

Before proving the extension theorem for $d<n$, we recall the definition of the Whitney decomposition of a set.
\begin{Def} A Whitney decomposition of an open set $A$ is a collection of closed cubes $\{Q_i\}_{i\in\N}$ with disjoint interiors and sides parallel to the axes such that $A=\bigcup_i Q_i$, where each cube $Q_i$ has side length $s_i=2^{-M_i}$ for some $M_i\in\mathbb Z$ and diameter $l_i$ such that $l_i\leq d(Q_i,A^c)\leq 4l_i.$
 \end{Def}
Denote by $x_i$ the center of the cube $Q_i$ and let $\varepsilon \in(0, \frac 1 4)$. Denote by $Q_i^* = (1 + \varepsilon)Q_i$ the cube with center $x_i$ expanded by factor $1+\varepsilon$. If $x\in Q_k\cap Q_i^*$ then  
\begin{equation}
1/4 s_k\leq s_i\leq 4s_k\label{trace:whitneysides}
\end{equation} 
and $Q_i$ and $Q_k$ touch each other. This implies that every point in $A$ is covered by at most $N_0$ cubes $Q_i^*,$ where $N_0\in\N$ depends only on $n$. By \cite[Section I.2.3]{jonwal} we can associate with decomposition $\{Q_i^*\}$ a partition of unity
$\{\varphi_i\}_{i\in\N}\subset C_c^\infty(\R^n)$, i.e.~a family of non-negative functions with the following properties:
\begin{align}
\begin{split}
&\supp\,\varphi_i\subset Q_i^*,\quad\sum \varphi_i=1\text{ on }\tilde D^c,\quad |D^j\varphi_i|\leq \tilde c l_i^{-|j|}\text{ for some }\tilde c>0.
\end{split}\label{trace:whitney}
\end{align}
If $D$ is an $d$-set in $\R^{n}$ then by \cite[Proposition VIII.1.1]{jonwal} the closure $\overline{D}$ of $D$ is also a $d$-set and $\mu(\overline{D}\setminus D)=0$ for every $d$-measure $\mu$. Therefore, it is enough to prove the theorem for a closed $d$-set $D$. 

Next, we define the extension operator $E$ from $H^{\psi,\alpha}(D,\mu)$ to $H^{\psi,\alpha}(\R^n)$, when $d<n$. Let $\omega_i=\mu(B(x_i,6l_i))^{-1}\asymp l_i^{-d}$ and $I=\{i\in\N:s_i\leq 1\}$. For $H^{\psi,\alpha}(D,\mu)$ define 
$$ E u(x)=\begin{cases} 
	u(x), & x\in  D \\        
    \displaystyle{\sum_{i\in
    I}\varphi_i(x)\omega_i\int_{|y-x_i|<6l_i}u(y)\mu(dy)},\quad & x\not\in  D.
    \end{cases}$$

\begin{Thm}\label{extension_theorem1} Let $D$ be a closed $d$-set in $\R^n$ and $d<n$. Let $\psi$ be a function such that \eqref{eq:j3}, \Hj and \Hd hold and $\alpha>0$. There exist a continuous extension operators $E:H^{\psi,\alpha}(D,\mu)\to H^{\psi,\alpha}(\R^n)$ such that $Eu=u$ $\mu$-a.e. on $D$ for all $u\in H^{\psi,\alpha}(D,\mu)$.
\end{Thm}
 \proof  
 We will show that for some $c>0$, $N\in\mathbb N$ and $k_0>\alpha\delta_2$ there exists a constant $\tilde c_1$ such that 
\begin{equation}
||\tilde Eu||_{(1),\alpha,k_0}^{h_0}\leq \tilde c_1||u||_{(2), D,\mu}^{c,N},\quad\forall u\in H^{\psi,\alpha}(D,\mu).\label{trace:show1}
\end{equation}
Here we can choose the smallest $k_0\in \N$ satisfying the given condition, see Theorem \ref{trace:moura} and Remark \ref{trace:Delta^k}. Since $D$ is of Lebesgue measure zero in $\R^{n}$ it is enough to prove \eqref{trace:show1} for $E u 1_{D^c}$. For every $x\in D^c$ there exists a $k$ such that $x\in Q_k$. If $s_k>4$ then by \eqref{trace:whitney} $x\not\in Q_i^*$ for all $i\in I$ and $Eu(x)=0$. Therefore it is enough to consider the case when $s_k\leq 4$. Also note that $\sum_i \varphi_i(x)=\sum_{i\in I}\varphi_i(x)$ when $s_k<1/4$.

Let $x\in Q_k$ and let $i\in I$ be such that $\phi_i(x)\neq 0$. Then for all $y\in B(x_i,6l_i)$ we have
\begin{equation}\label{ball29}
|y-x_k|\leq|y-x_i|+|x_i-x|+|x-x_k|\leq 6l_i+l_i+l_k\leq 29 l_k,
\end{equation}
which implies that
\begin{align}
|Eu(x)|&\leq\sum_{i\in I}\varphi_i(x)\omega_i\int_{|y-x_i|<6l_i}|u(y)|\mu(dy)\lesssim\sum_{i\in I}\varphi_i(x)l_i^{-d}\int_{|y-x_k|<29l_k}|u(y)|\mu(dy)\nonumber\\
&\overset{\eqref{trace:whitneysides}}{\lesssim} \sum_{i\in I}\varphi_i(x) \cdot
l_k^{-d}\int_{|y-x_k|<29l_k}|u(y)|\mu(dy)\lesssim \left(
l_k^{-d}\int_{|y-x_k|<29l_k}u^2(y)\mu(dy)\right)^{1/2}\nonumber.
\end{align}
For $j\in\N$ define $\Delta_j:=\bigcup_{\{k:s_k=2^{-j}\}}Q_k$. Note that there exists an integer $N_1$ depending only on $n$ such that every point $y\in \R^n$ is covered by at most $N_1$ balls $B(x_k,29l_k)$ where $Q_k\subset \Delta_j$. This follows from the fact
that $|x_k-x_{k'}|\geq 2^{-j}$ and $l_k=\sqrt{n} 2^{-j}$, for all
$Q_k,Q_{k'}\subset \Delta_j$. By the previous calculation it follows that 
\begin{align*}
\int_{ D^c}| Eu(x)|^2dx&=\sum_{j=-2}^\infty\sum_{Q_k\subset\Delta_j}\int_{Q_k}| Eu(x)|^2dx\lesssim\sum_{j=-2}^\infty\sum_{Q_k\subset \Delta_j}\int_{Q_k} \left( l_k^{-d}\int_{|y-x_k|<29l_k}u^2(y)\mu(dy)\right)dx\nonumber\\
& \asymp\sum_{j=-2}^\infty 2^{-(n-d)j} \sum_{Q_k\subset \Delta_j}\int_{|y-x_k|<29l_k}u^2(y)\mu(dy) \leq 2^{2(n-d)}\sum_{j=-2}^\infty \int_{\Delta_j}u^2(y)\mu(dy)\nonumber,
\end{align*} 
which implies that
\begin{equation}
 ||Eu||_{L^2(\R^{n})}\lesssim ||u||_{L^2( D,\mu)}. \label{trace:proof1}
\end{equation}
Next, for $x\in \Delta_{i}$, $y\in\Delta_j$ and
$|x-y|<2^{-i}/2$ we have
\begin{align}
2^{-j}\sqrt{n}&\leq d(\Delta_j, D)\leq d(y,D)\leq |x-y|+d(x, D)\leq 6\sqrt{n} 2^{-i},\nonumber
\end{align}
so $j\geq i-2$. Analogously, $\sqrt{n} 2^{-i}\leq 5\sqrt{n} 2^{-j}+2^{-i-1}$ so $j\leq i+2$. Therefore,  
\begin{equation}
\label{trace:proof2}x\in \Delta_{i},\,|x-y|<2^{-i}/2\,\Rightarrow\, y\in\bigcup_{j=i-2}^{i+2}\Delta_j.
\end{equation}
Since $Eu(x)=0$ for for $x\in\Delta_i$, $i\leq -3$ it follows that $Eu(y)=0$ when $|x-y|<2^{-i-1}$ for some $x\in\Delta_i$, $i\leq -5$. Analogously, $\Delta_h^{k_0} (\tilde Eu)(x)=0$ if $|h|<2^{-5}/k_0$. Therefore, for $h_i:=2^{-i-1}/k_0$,
\begin{align}
\int\limits_{D^c}\int\limits_{|h|<h_0}|&\Delta_h^{k_0} (Eu)(x)|^2\frac{\psi^\alpha\left(|h|^{-1}\right)}{|h|^{n}}
dh\,dx\nonumber\leq
\sum_{i=-4}^\infty\iint\limits_{\substack{x\in\Delta_i\\|h|<h_i}}|\Delta_h^{k_0} (Eu)(x)|^2\frac{\psi^\alpha\left(|h|^{-1}\right)}{|h|^{n}}
dh\,dx\nonumber\\
&+\sum_{i=5}^\infty\iint\limits_{\substack{x\in\Delta_i\\h_{i}\leq|h|<2^{-5}/k_0}}|\Delta_h^{k_0} ( Eu)(x)|^2\frac{\psi^\alpha\left(|h|^{-1}\right)}{|h|^{n}}
dh\,dx=:A+B\nonumber.
\end{align}
First we asses the term $B$. Let $F_{i}:=\bigcup_{j=i}^\infty \Delta_j.$ Note that
\begin{align} 
B&=\sum_{i=5}^\infty\sum_{m=4}^{i-1}\int\limits_{h_{m+1}\leq|h|<h_{m}}\frac{\psi^\alpha\left(|h|^{-1}\right)}{|h|^{n}}\int\limits_{x\in\Delta_i}|\Delta_h^{k_0} (Eu)(x)|^2
dx\,dh\nonumber\\
&=\sum_{m=4}^{\infty}\int\limits_{h_{m+1}\leq|h|<h_{m}}\frac{\psi^{\alpha}\left(|h|^{-1}\right)}{|h|^{n}}\int\limits_{x\in F_{m+1}}|\Delta_h^{k_0} ( Eu)(x)|^2dx\,dh\nonumber\\
&\overset{\eqref{eq:phi}}{\lesssim} 
\sum_{m=4}^{\infty}\psi^\alpha\left(2^{m}\right)2^{mn}\iint\limits_{\substack{x\in F_{m+1}\\h_{m+1}\leq|h|<h_{m}}} |\Delta_h^{k_0} (Eu)(x)|^2 dx\,dh\nonumber.
\end{align}
Similarly as in \eqref{trace:proof2}, for $x\in F_{i+1}$ and $|h|<h_i$ it follows that ${x,x+h,..,x+k_0h\in F_{i-2}}$. Since 
$$|\Delta_h^{k_0} (Eu)(x)|^2\lesssim | Eu(x)- Eu(x+h)|^2+\cdots +|Eu(x+(k_0-1)h)-Eu(x+k_0h)|^2$$
it follows that
\begin{align}\label{trace:proofB}
B&\lesssim \sum_{m=4}^{\infty}\psi^\alpha\left(2^{m}\right)2^{mn}
\iint\limits_{\substack{x,y\in F_{m-2}\\|x-y|<h_m}}(Eu(x)-Eu(y))^2 dx\,dy.
\end{align}
For $k,m\geq 2$ and $x\in\Delta_k$ and $y\in\Delta_m$ it follows that 
\begin{align*}
| Eu(x)-Eu(y)|\leq\sum_i\sum_j\varphi_i(x)\varphi_j(y)\omega_i\omega_j\iint\limits_{\substack{|s-x_i|<6l_i\\|t-x_j|<6l_j}}|u(s)-u(t)|\mu(ds)\mu(dt)
\end{align*}
From
\begin{equation}\label{eq:diam_comparable}
x\in\Delta_k,\,\varphi_i(x)\neq 0\Rightarrow \frac 1 8 l_k\leq l_i\leq 64 l_k
\end{equation}
it follows that 
\begin{align*}
|Eu(x)-Eu(y)|&\overset{\eqref{ball29}}{\lesssim}\sum_i\sum_j\varphi_i(x)\varphi_j(y)l_k^{-d}l_m^{-d}\iint\limits_{\substack{|s-x_k|<29l_k\\|t-x_m|<29l_m}}|u(s)-u(t)|\mu(ds)\mu(dt)\\
&\leq \left(l_k^{-d}l_m^{-d}\iint\limits_{\substack{|s-x_k|<29l_k\\|t-x_m|<29l_m}}(u(t)-u(s))^2\mu(dt)\mu(ds)
\right)^{1/2}.
\end{align*}
Here $x_k,x_m$ are the centers and $l_k,l_m$ diameters of cubes $Q_p\subset\Delta_k$ and $Q_r\subset\Delta_m$ containing $x$ and $y$ respectively. 
Now it follows that for $i\in\N$, $y\in\Delta_m$ and $k,m\geq 2$ 
\begin{align*}
\int_{\substack{x\in\Delta_k,\\|x-y|<2^{-i}}}|Eu(x)-Eu(y)|^2 dx&\lesssim
\int_{\substack{x\in\Delta_k\\|x-y|<2^{-i}}}l_k^{-d}l_m^{-d}\iint\limits_{\substack{|s-x_k|<29l_k\\|t-x_m|<29l_m}}(u(t)-u(s))^2\mu(dt)\mu(ds)dx\\
&\leq N_0 l_k^{-d}l_m^{-d}s_k^{-n}\iint\limits_{\substack{|s-y|<c2^{-k}+2^{-i}\\|t-x_m|<29l_m}}(u(t)-u(s))^2\mu(dt)\mu(ds),
\end{align*}
where $c=30\sqrt{n}$. Analogously, we get
$$\iint\limits_{\substack{x\in\Delta_k,\,y\in\Delta_m\\|x-y|<2^{-i}}}|Eu(x)-Eu(y)|^2
dx dy \lesssim N_0^2 2^{-(n-d)k}2^{-(n-d)m}\iint\limits_{\mathclap{|t-s|<2^{-i}+c2^{-k}+c2^{-m}}}(u(t)-u(s))^2d\mu(s)d\mu(t).$$
This implies that for $i\geq 4$ 
\begin{align}
&\iint\limits_{\substack{x,y\in F_{i-2}\\|x-y|<2^{-i}}}| Eu(x)- Eu(y)|^2
dx\,dy=\sum_{k,m=i-2}^\infty\,\, \iint\limits_{\substack{x\in\Delta_k,\,y\in\Delta_m\\|x-y|<2^{-i}}}| Eu(x)- Eu(y)|^2
dx dy\nonumber\\
&\lesssim\left(\sum_{k,m=i-2}^\infty
2^{(d-n)k}2^{(d-n)m}\right)\iint\limits_{|t-s|<(8c+1)2^{-i}}(u(t)-u(s))^2\mu(ds)\mu(dt)\nonumber
\end{align}
and by applying this to \eqref{trace:proofB} we get that
\begin{align}
B& \lesssim \sum_{i=4}^\infty \psi^\alpha\left(2^{i}\right) 2^{i(2d-n)}
\iint\limits_{\mathclap{|t-s|<(8c+1)2^{-i}}}(u(t)-u(s))^2\mu(ds)\mu(dt).\label{trace:proof3}
\end{align}
Next, by the mean value theorem 
\begin{align}
A&\,\lesssim\sum_{i=-4}^\infty\,\int\limits_{|h|<h_{i}}\left(\sum_{|j|=k_0}\,\int\limits_{\Delta_i}\int_0^1...\int_0^1
|h|^{2k_0} |D^{j}(Eu)(x+(t_1+...+t_{k_0})h)|^2 dt_1\,...dt_{k_0}\, dx\right)\frac{\psi^\alpha\left(|h|^{-1}\right)}{|h|^{n}}dh\nonumber\\
&\overset{\eqref{trace:proof2}}{\leq} \sum_{i=-4}^\infty\,\int\limits_{|h|<h_{i}}\frac{\psi^\alpha\left(|h|^{-1}\right)}{|h|^{n-k_0}}dh\cdot\sum_{|j|=k_0}\int\limits_{F_{i-2}\setminus F_{i+3}}
|D^{j}(Eu)(z)|^2 dz\overset{\Hj}{\lesssim }5\sum_{i=-2}^\infty \frac{\psi^\alpha(2^i)}{2^{2k_0i}}\sum_{|j|=k_0}\int_{\Delta_i} |D^j(Eu)(z)|^2 dz\label{trace:proofA}.
\end{align}
In the last line we also used that $D^{j}(\tilde Eu)(z)=0$ if $z\in\Delta_i$ and $i\leq -3$. To find an upper bound for $|D^j(Eu)(z)|$, $z\in\Delta_i$, we distinguish two cases; $i\geq 2$ and $i<2$. First, take $z,y\in Q_m\subset\Delta_l$, $l\geq 2$ and $|j|=k_0$. Since $\displaystyle{\sum_iD^j\varphi_i(z)=0}$, by similar calculations as before, we get
\begin{align*}
|D^j( Eu)(z)|&=\left|\sum_iD^j\varphi_i(z)\omega_i\int\limits_{\mathclap{|s-x_i|<6l_i}}(u(s)-\tilde Eu(y))\mu(ds)\right|\\
&\leq\sum_i\sum_k |D^j\varphi_i(z)|\varphi_k(y)\left(\omega_i\omega_k\iint\limits_{\mathclap{\substack{|s-x_i|<6l_i\\|t-x_k|<6l_k}}}|u(s)-u(t)|^2\mu(ds)\mu(dt)\right)^{\frac 1 2}.
\end{align*}
Recall that there are at most $N_0$ indices $i$ for which $z\in Q_i^*$ and $D^j\varphi_i(z)\neq 0$. By \eqref{trace:whitneysides} and \eqref{trace:whitney} $z\in Q_i^*$ implies $\omega_i\asymp l_m^{-d}$ and $|D^j\varphi_i(z)|\lesssim l_i^{-|j|}\lesssim  l_m^{-|j|}$. Also, by \eqref{eq:diam_comparable} $\omega_k\asymp l_m^{-d}$ for $k$ such that $\varphi_k(y)\neq 0$. Therefore, 
\begin{align*}
|D^j( Eu)(z)|&\overset{\eqref{ball29}}{\lesssim}l_m^{-k_0}\left(l_m^{-2d}\iint\limits_{\mathclap{\substack{|s-x_m|<29l_m\\|t-x_m|<29l_m}}}|u(s)-u(t)|^2\mu(ds)\mu(dt)\right)^{\frac 1 2}.
\end{align*}
Applying this inequality to \eqref{trace:proofA} we arrive to 
\begin{align}
A&\lesssim \sum_{i=2}^\infty \psi^\alpha(2^i)2^{-2k_0i}\sum_{|j|=k_0}\sum_{Q_m\subset \Delta_i}\int_{Q_m}  2^{2k_0i+2id}\iint\limits_{\mathclap{\substack{|s-x_m|<29l_m\\|t-x_m|<29l_m}}}|u(s)-u(t)|^2\mu(ds)\mu(dt)dz\nonumber
\end{align}
and since every $s\in D^c$ is covered by at most $N_1$ balls $B(x_m,29l_m)$ it follows that
\begin{align}
A\lesssim&\sum_{i=2}^\infty \psi^\alpha(2^i)2^{i(2d-n)}\iint\limits_{\mathclap{|s-t|<60\sqrt{n}2^{-i}}}|u(s)-u(t)|^2\mu(ds)\mu(dt). \label{trace:proof4} 
\end{align}
For the remaining part in $A$, take $z\in\Delta_k$, $k\geq -2$. By the same arguments as before, 
\begin{align*}
|D^j( Eu)(z)|&\leq\sum_i|D^j\varphi_i(z)|\omega_i\int\limits_{\mathclap{|s-x_i|<6l_i}}|u(s)|\mu(ds)\overset{\eqref{trace:whitney}}{\lesssim} \sum\limits_{\varphi_i(z)\neq 0}l_i^{-2}\left(\omega_i\int\limits_{\mathclap{|s-x_i|<6l_i}}|u(s)|^2\mu(ds)\right)^{\frac 1 2}\\
&\lesssim 2^{kk_0}\left(2^{kd}\int\limits_{\mathclap{|s-x_k|<29l_k}}|u(s)|^2\mu(ds)\right)^{\frac 1 2}
\end{align*}
and therefore $\sum_{i=-2}^1\int_{\Delta_i}|D(Eu)(z)|^2 dz\lesssim||u||_{L^2( D,\mu)}$. This inequality together with \eqref{trace:proof1}, \eqref{trace:proof3} and \eqref{trace:proof4} implies \eqref{trace:show1}.

That $E$ is truly the extension operator for $R$, i.e. that $REu =u$ $\mu$-a.e., follows from calculation similar to the calculation above. One first shows that for every $t_0\in D$ and $r>0$ small enough
\begin{align}
&\int_{|x-t_0|\leq r}(Eu(x)-u(t_0))^2dx \lesssim \frac{r^{d}}{\psi^\alpha(r^{-1})} \int_{|t-t_0|<30 r}(u(t)-u(t_0))^2\frac{\psi^\alpha(|t-t_0|^{-1})}{|t-t_0|^{2d-n}} \mu(dt)\nonumber,
\end{align}
where the integral is finite for $\mu$-almost all $t_0$ and decreasing as $r$ goes to 0. Since $\displaystyle{\lim\limits_{r\to 0}\frac{1}{\psi^\alpha(r^{-1})r^{n-d}}=0}$, it follows that for $\mu$ almost all $t_0$
\begin{align}
&|R Eu(t_0)-u(t_0)|\lesssim \lim_{r\to 0} \left(r^{-n}\int_{|x-t_0|\leq 
r}(\tilde Eu(x)-u(t_0))^2dx\right)^{1/2}= 0\nonumber.
\end{align}
\qed

\begin{Rem}\label{trace_rem}
Let $\alpha\in(\frac 1 2,\infty)$. Note that for a $n$-set $D$ in $\R^n$, $\tilde D=D\times\{0\}$ is a $n$-set in $\R^{n+1}$ and that every function $u\in H^{\psi,\alpha}(D,\mu)$ can be represented as a function $\tilde u$ in $H^{\tilde\psi,\tilde\alpha}(\tilde D,\tilde\mu)$, where
\begin{align*}
&\tilde\alpha:=2\alpha,\ \tilde\psi(|\xi|):=\psi^\frac 1 2(|\xi|)|\xi|^\frac1{\tilde\alpha},\ \tilde \mu(A\times\{0\}):=\mu(A)\\
&\tilde u (x,0):=u(x),\, x\in D\quad
\text{and}\quad ||\tilde u||_{(1),\tilde D,\tilde\mu}=|| u||_{(1),D,\mu}, 
\end{align*}
By \cite[Theorem 7.13.]{bernstein} the function $\xi\mapsto\phi^\frac 1 2 (|\xi|)|\xi|^\frac 1{2\alpha}$ is a complete Bernstein function and $\tilde\psi$ satisfies conditions \Hj and \Hd with $\tilde\delta_i:=\frac{\delta_i}2+\frac{1}{4\alpha}$. Analogously, the space $H^{\psi,\alpha}(\R^n)$ can be represented as $H^{\tilde\psi,\tilde\alpha}(\R^{n}\times\{0\},\bar\mu)$, where $\bar\mu$ is the restriction of the $n$-dimensional Hausdorff measure in $\R^{n+1}$ to $\R^n\times\{0\}.$\\
\end{Rem}

\begin{Thm}\label{extension_theorem3}
Theorem \ref{extension_theorem1} holds true as well in the case of $d=n$ and $\alpha\in(\frac 1 2,\infty)$. 
\end{Thm}

\proof
Take $u\in H^{\psi,\alpha}(D,\mu)$ and let $\tilde u$ be the corresponding function in $H^{\tilde\psi,\tilde\alpha}(\tilde D,\tilde\mu)$ from Remark \ref{trace_rem}. By Theorem \ref{extension_theorem1} function $\tilde u$ can be extended to a function $\tilde{E} \tilde u\in H^{\tilde\psi,\tilde\alpha}(\R^{n+1})$, which can then be restricted to a function in $H^{\tilde\psi,\tilde\alpha}(\R^n\times\{0\},\bar\mu)$ by applying the continuous restriction operator $\bar R$ from Theorem \ref{restriction_theorem}. Again using the Remark \ref{trace_rem}, we can define the extension operator $E$ as
$$(Eu)(x)=(\bar R\tilde{E} \tilde u)(x,0),\,x\in\R^n.$$
The continuity of $E$ follows from the continuity of the extension and restriction operators $\tilde{E}$ and $\bar R$ and $REu=u$ almost everywhere on $D$, where $R$ is the restriction operator from $H^{\psi,\alpha}(\R^n)$ to $H^{\psi,\alpha}(D,\mu)$.
\qed

\proofof{{\bf Corollary \ref{corollary}}:}
By applying Theorem \ref{thm:main}, the proof is an immediate consequence of \cite[Theorem 1.1, Corollary 2.9]{wagner2}.
\qed

\bibliographystyle{plain}
\bibliography{tracethm}

\end{document}